\documentclass[11pt,twoside]{article}
\usepackage{latexsym,amssymb,amsfonts}
\usepackage{amsmath}
\usepackage{amscd}
\usepackage[english]{babel}
\usepackage{pstricks}
\usepackage{pst-node}
\usepackage{multido}
\usepackage{pstricks-add}
\usepackage{pst-math}
\begin{document}

\title{Frobenius algebras and skein modules \\ 
of surfaces in $3$-manifolds}
\author{Uwe Kaiser \footnote{partially
supported by NSF grant DMS-0204627}\\ 
  \begin{tabular}{c}
    Boise State University\\
  \texttt{kaiser@math.boisestate.edu}
\end{tabular}}

\date{February 27, 2008}

\maketitle

\abstract{\noindent For each Frobenius algebra there is defined a skein 
module of surfaces embedded in a given $3$- manifold and bounding a  
prescribed curve system in the boundary.
The skein relations are local and generate 
the kernel of a certain natural extension of the corresponding 
topological quantum field theory.
In particular the skein module of the $3$-ball is isomorphic to the 
ground ring of the Frobenius algebra.
We prove a presentation theorem for the skein module with 
generators incompressible surfaces colored by elements of a generating set of 
the Frobenius algebra, and with relations determined by tubing geometry 
in the manifold and relations of the algebra.}

\vskip .2in


\vskip .1in

\noindent \textsf{Mathsubject Class.:} 57M25, 57M35, 57R42

\vskip .3in

\section{Introduction.}

In \cite{P} the idea of algebraic topology based on knots originated.
The idea is to study quotients of formal linear combinations of links
by local relations that allow computation of the modules and understanding of the interplay between link theory and geometric topology of the $3$-manifold. 
 Most interesting skein modules of links in $3$-manifolds are based on skein 
relations with the property that the modules of the $3$-ball are 
isomorphic to the ground ring. Typically one considers relations 
suggested by link polynomials like the Jones, Homflypt or Kauffman 
bracket polynomials. Understanding the structure of these modules for 
general $3$-manifolds is an important question which is related to 
basic problems of quantum topology like the volume conjecture
\cite{L}. Interestingly, these polynomial link invariants invariants have 
operator-valued
extensions as $1$-dimensional embedded topological quantum 
field theories \cite{T1}, \cite{T2}.
  
\vskip .05in

It is well-known that theories of surfaces in $3$-manifolds provide
  powerful tools towards the understanding
 of $3$-manifold topology. So it seems to be almost surprising that 
general skein theory of surfaces in $3$-manifolds have not yet been developed.
Probably the reason is that embeddings of surfaces appear to be much more 
rigid. More importantly, the \textit{local theory} seems to be of 
no interest because there are no nontrivial \textit{local} embeddings of 
surfaces 
while knot theory is interesting because of the local possibility of 
knotting and linking.

\vskip .05in

Recently the study of skein relations between surfaces was suggested 
by Bar-Natan \cite{BN}. It has been put forward into the framework of 
general $3$-manifolds by Asaeda and Frohman \cite{AF}. The structure
of skein modules and related categories have been studies by Gad Naot
\cite{N1}, \cite{N2}. In \cite{K1} Khovanov initiates the discussion of skein relations
between surfaces originating from more general Frobenius algebras.
The abstract $2$-dimensional 
topological quantum field theories defined from Frobenius algebras
turn out to be mostly interesting for Khovanov theory, and their 
role in general $3$-manifold theory remains to be explored, see \cite{FK}. 
But we will see that in general \textit{globalizing} the corresponding
skein relations of surfaces yields structures which not surprisingly
measure \textit{global} embedding phenomena of surfaces in $3$-manifolds
in a nontrivial way.    

\vskip .05in

Until now, because of the emphasis on the Khovanov theory background, many 
results have been proven separately for different skein relations. 
It is the goal of this paper to initiate a discussion of a large class 
of skein modules defined from embedded surfaces in $3$-manifolds.
We will see that many properties of modules
defined using Frobenius algebras are closely related to the corresponding $2$-dimensional topological quantum field theories. 
In fact we will show that the skein relations naturally emerge from a 
discussion of the kernel of a suitably \textit{extended} topological quantum 
field theory. Also, many results proven for 
specific skein relations will hold in broad generality. In section 9 we 
prove a presentation theorem which shows that in the case of free 
Frobenius algebra the module can be generated by surfaces colored with elements of the basis, and relations defined from tubing geometry of the $3$-manifold. 
For example the tunnel number of a knot naturally appears in the discussion 
of skein modules of surfaces in knot complements. 
In general, the relations of our skein modules of surfaces in $3$-manifolds with respect to natural generating sets, can be expressed through the difficult geometric problem of understanding the different ways to compress a given 
surface in a $3$-manifold.

\vskip .05in

Frobenius algebras also naturally appear in the study of quantum cohomology, and
through the study of Frobenius manifolds are related with loop spaces, 
symplectic geometry and Landau-Ginzburg models, see \cite{H} for a good introduction. 

\vskip .05in

Most of the ideas that we discuss here are already present in the work of Asaeda and Frohman \cite{AF}.
But we will elaborate on various subtleties emerging in the discussion 
of skein modules of surfaces in $3$-manifolds from general Frobenius algebras,
which do not seem to appear elsewhere in the published literature.
The author is aware that most of the material in this paper is folklore in the community of Khovanov theory researchers. 

\vskip .05in

I would like to thank Masahico Saito for suggesting the study of skein modules of surfaces for general Frobenius algebras. I also would like to thank Carmen Caprau and Charlie Frohman for several interesting and helpful remarks.
It has been Charlie's enthusiasm convincing me to proceed with this subject.

\section{Frobenius algebras and Frobenius systems.}

For further details about some of the following algebraic concepts 
see \cite{K1} and \cite{K}. We like to point out what we call a Frobenius algebra is called a Frobenius system in \cite{K1}.    

Throughout let $R$ be a commutative ring with $1$. All tensor products are 
over $R$, so $\otimes =\otimes_R$.

A \textit{Frobenius algebra} 
$\mathcal{F}=(R, A, \Delta, \varepsilon)$ is a commutative ring $A$ with 
$1$ and inclusion $\iota : R\rightarrow A$ such that $\iota (1)=1$.
Thus $A$ is an $R$-algebra with product
$$\mu : A\otimes A\rightarrow A$$
and unit $\iota $. Also
$$\Delta : A\rightarrow A\otimes A$$
is a cocommutaive and coassociative coproduct, which is 
an $A$-bimodule map, and has counit 
the $R$-module map $\varepsilon : A\rightarrow R$.   

Note that the counit condition means
$$(\varepsilon \otimes Id)\circ \Delta =Id,\leqno(1)$$
where we identify $R\otimes A=A\otimes R=A$ using the restrictions 
of $\mu $. 
 
The $R$-algebra $A$ is called a \textit{Frobenius extension} of 
$R$.

If we write $\Delta (1)=\sum_{i=1}^ru_i\otimes v_i$ for 
$u_i,v_i\in A$ then $(1)$ and cocommutativity imply that for all $a\in A$: 
$$a=\sum_{i=1}^r\varepsilon (au_i)v_i=\sum_{i=1}^r\varepsilon (av_i)u_i,\leqno(2)$$
In particular, $A$ is a finitely generated $R$-module with the two generating sets 
$\{u_i|i=1,\ldots r\}$ and $\{v_i|i=1,\ldots ,r\}$.
 
Conversely, each set $\{(u_i,v_i)|i=1,\ldots r\}\subset A\times A$
such that $\sum_{i=1}^ru_i\otimes v_i=\sum_{i=1}^rv_i\otimes u_i$ and 
$(2)$ holds for a given $\varepsilon $, defines a cocommutative 
and coassociative $A$-bimodule coproduct by
$$\Delta (1)=\sum_{i=1}^ru_i\otimes v_i$$
with counit $\varepsilon $. 
 
Following \cite{K} we call $(R,A,\{(u_i,v_i)|i=1,\ldots ,r\}, \varepsilon )$ a
\textit{Frobenius system}. Each Frobenius system determines a 
Frobenius algebra, and conversely for each Frobenius algebra there 
exists a not necessarily unique Frobenius system.
The \textit{rank} of a Frobenius system 
is defined by $\sum_{i=1}^r\varepsilon (u_iv_i)=\varepsilon \mu \Delta (1)\in R$,
which only depends on the underlying Frobenius algebra.
 
\vskip .1in

\noindent \textbf{Remark 2.1.} (a) Let $\{(u_i,v_i)|i=1,\ldots r \}$ be defining a Frobenius 
system on the $R$-algebra $A$ and let $\{w_i|i=1,\ldots ,s\}$ be a generating set of $A$. Then there exist $\lambda_{ij}\in R$ for
$i=1,\ldots ,r$ and $j=1,\ldots ,s$ such that $v_i=\sum_{j=1}^s\lambda_{ij}w_j$. Thus
$\{(z_j,w_j)|j=1,\ldots s\}$ with $z_j=\sum_{i=1}^r\lambda_{ij}u_i$ define a usually different Frobenius systems with the same underlying Frobenius algebra. 
 
\noindent (b) Usually there is the additional assumption that $A$ is a 
projective $R$-module, i.\ e.\ summand of a free $R$-module. This condition plays a role in discussions of the structure theory of topological quantum field theories 
associated to Frobenius algebras. It will not at all be important in our discussions.

\noindent (c) To avoid redundancies one can assume  
in the definition of Frobenius systems that for each subset 
$J\subset \{1,\ldots ,r\}$
$$\sum_{i\in J}u_i\otimes v_i\neq 0.$$

\vskip .1in

The natural $R$-module map
$$A\ni a\mapsto (A\ni b\mapsto \varepsilon (ba)) \leqno(3)$$
defines an $R$-isomorphism 
$$A\rightarrow A^*=Hom_R(A,R).$$

If $A$ is a free $R$-module then we can find a Frobenius system 
for it with both $\{u_i\}$ and $\{v_i\}$ $R$-bases.
Then by (3) applied to the functionals mapping some $u_i$ to $1$ and all 
other $u_j$ to $0$. 
it follows that there exists a \textit{dual} basis $\{w_i\}$ of $A$ such 
that 
$$\varepsilon (u_iw_j)=\delta_{ij},$$
where $\delta_{ij}=1$ for $i=j$ and $=0$ for $i\neq j$. 
A \textit{free Frobenius algebra} is a Frobenius algebra with $A$ a 
free $R$-module.

\vskip .1in

In Khovanov theory gradings are important. A Frobenius algebra is \textit{graded} if $A$ is a graded ring (a graded abelian group with $\mu $ a degree $0$ mapping) such that all structure maps are graded.

\vskip .1in 

\noindent \textbf{Examples 2.2.} (a) Following \cite{K1} we define 
the \textit{universal rank two Frobenius algebra} with $R=\mathbb{Z}[h,t]$, $A=R[X]/(X^2-hX-t)$ and $\varepsilon (1)=0$, $\varepsilon (X)=1$
and \textit{structure elements} \\ $\{(1,X-h),(X,1)\}$
corresponding to the bases $\{1,X\}$ and $\{X-h,1\}$.
Then
$$\Delta (1)=1\otimes X+X\otimes 1-h1\otimes 1.$$
It follows from $A$-bilinearity that:

\vskip -.3in

\begin{align*}
\Delta (X)=\Delta (1)X&=(1\otimes X+X\otimes 1-h(1\otimes 1))X\\
&=1\otimes X^2+X\otimes X-h(1\otimes X)\\
&=1\otimes (hX+t)+X\otimes X-h(1\otimes X)\\
&=X\otimes X+t(1\otimes 1)
\end{align*}

\vskip .1in

 For $a\in A$ let $\mu_a: A\rightarrow A$ be defined by
$$\mu_a(c)=\mu (a\otimes c)=ac$$
for $c\in A$. 

\vskip .1in

Note that the values $\varepsilon_n:=\varepsilon (X^n)\in R$
are recursively determined by $\varepsilon_0=0, \varepsilon_1=1$ and
$$\varepsilon_{n+2}=h\varepsilon_{n+1}+t\varepsilon _n$$
for all $n\geq 0$. The corresponding universal Frobenius algebra
$\mathcal{F}_{\mathbf{U}}$ is graded by the $deg(X)=1$ 
and $deg(h)=2$, $deg(t)=4$.

\noindent (b) Let $R=\mathbb{Z}$ and $G$ a finite abelian group. Then 
the group algebra $\mathbb{Z}[G]$ is a Frobenius extension of $\mathbb{Z}$. 
The natural Frobenius algebra $\mathcal{F}_G$ is defined by 
$\Delta (1)=\sum_gg^{-1}\otimes g$, with  
$\varepsilon (g)=0$ for $g\neq e$ and $\varepsilon (e)=1$,
where $e\in G$ is the identity of the group.
This is a free Frobenius algebra with basis $G$.

\noindent (c) Given a Frobenius algebra $\mathcal{F}=(R,A,\Delta ,\varepsilon)$
and some invertible element $y\in A$ there is defined the \textit{twisted}
Frobenius algebra $\mathcal{F}^y:=(R,A,\Delta ',\varepsilon ')$ by
$\Delta '=\Delta \circ \mu_{y^{-1}}$ and 
$\varepsilon '=\varepsilon \circ \mu_y$. It is known that any two 
Frobenius algebras defining the same Frobenius extension differ by twisting,
see \cite{K}, Theorem 1.6. 

\vskip .1in

Given a Frobenius algebra $\mathcal{F}$ we let $M(\mathcal{F})$ denote
the symmetric monoidal subcategory of the category of all (projective) 
$R$-modules with objects $A^{\otimes r}, r\geq 0$ ($A^0:=R$), 
which is generated by tensor product and composition from 
the \textit{structure morphisms} $\mu , \iota , \Delta , \varepsilon ,\tau $
and $\mu _a$ for all $a\in A$.  
Note that all \textit{structure morphisms} are $R$-module maps and only $\iota $ and 
$\varepsilon $ are not $A$-module respectively $A$-bimodule maps.
The set of objects is naturally identified with the set of natural 
numbers $\{0,1,2,\ldots \}$.

\section{Extending $2$-dimensional topological quantum field \\ theories.}

It is well-known that Frobenius algebras (with ground ring $R$) define (orientable) $2$-dimensional 
topological quantum field theories 
with values in the category of $R$-modules.
Let $Cob$ denote the category with objects $1$-dimensional closed manifolds, standardly embedded in $\mathbb{R}^2$ 
(thus the objects are naturally identified with the natural numbers).
The morphisms are 
isotopy classes of orientable $2$-manifolds, 
embedded in $\mathbb{R}^N\times [0,1]$, for $N$ large,
bounding corresponding input respectively output $1$-manifolds
in $\mathbb{R}^2\times 1$ respectively $\mathbb{R}^2\times 0$.
The monoidal structure is defined by disjoint unions and the symmetric
structure is defined by the switch surface defined from two cylinders.

For $r\geq 0$ the topological quantum field theory assigns to the $1$-manifold $\amalg_rS^1$ ($\amalg $ is disjoint union) the $R$-module
$A^{\otimes r}$ with $A^{\otimes 0}:=R$.  
 The the structure maps $\Delta, \mu, \varepsilon, \iota$ and $\tau$ 
of the Frobenius algebra are assigned to the 
the pair of pants \textit{as you wear it}, the pair of pants \textit{turned upside down}, the cap (local maximum) , and the cap upside down (local minimum)
and the disjoint union of two cylinders switching the order of two input and output circles.

This defines the symmetric momoidal functor:
$$F(\mathcal{F}): Cob \rightarrow M(\mathcal{F})$$
with values in the symmetric monoidal subcategory of $M(\mathcal{F})$, 
which is monoidally generated by the structure maps of $\mathcal{F}$
(in fact in the subcategory of $\mathcal{M}(\mathcal{F})$ generated by the structure morphisms except the $\mu_a$).

\vskip .1in

Recall the usual way to calculate the $R$-module morphism assigned to a
surface by $F(\mathcal{F})$: Approximate the height function 
$S\subset \mathbb{R}^N\times [0,1]\rightarrow [0,1]$ by 
a Morse function 
with values in $[0,1]$
We assume that it takes
value $1$ on input circles and value $0$ on output circles. 
Then break up $[0,1]$ into subintervals, each containing 
only one critical point, the $R$-module map is determined 
by the monoidal functor property of the topological quantum field theory 
from the values on the standard surfaces. 
Finally represent possible permutations by sequences of transpositions
of neighboring circles.
This is in general a 
difficult calculation if carried out explicitly. 

It is essentially Bar-Natan's idea \cite{BN}, at least in a special case, to realize that the calculation can be carried out \textit{locally} using skein relations. This is important if one tries to  understand what properties are necessary such that Khovanov homology can be defined from a particular Frobenius algebra, see \cite{K1}. 

We will extend his idea to the general case of Frobenius algebras. 

\vskip .1in

\noindent \textbf{Definition 3.1.} Let $A$ be a commutative ring with $1$. Then we define the category
$Cob(A)$ as follows: The objects are the same as those of $Cob$. The morphisms are isotopy classes of surfaces with components colored by elements of $A$.
Using the multiplication of $A$ the composition of morphisms is defined.
Moreover $Cob\subset Cob(A)$ is the obvious subcategory of those morphisms 
with only colors $1$.

\vskip .1in

The morphisms of $Cob(A)$ can also be considered as pairs $(S,w)$ where $S$ is an object of $Cob$ and $w\in H^0(S,A)=Hom(H_0(S),A)$ is the coloring.

\vskip .1in

\noindent \textbf{Theorem 3.2.} The functor $\mathcal{F}$ naturally extends to a 
symmetric monoidal functor, also denoted 
$$F(\mathcal{F}): Cob(A)\rightarrow M(\mathcal{F}).$$

\noindent \textit{Proof.} 
Using monoidal structures it suffices to define the 
fucntor for connected surfaces in $Cob(A)$. 
Suppose a connected orientable surface $S$ is colored by $a\in A$.
First assume that the boundary of $S$ is not empty.
Then choose an input or output circle and 
pre- or postcompose by the corresponding map $Id \otimes \ldots \otimes \mu_a \otimes \ldots \otimes Id$. We know that $\Delta $ respectively $\mu $ are $A$-bimodule maps and commutative respectively cocommutative.
Thus it does not matter which input or output circle is chosen. In fact, if $a=a_1a_2\ldots a_{r+s}$ is any factorization then the factors can be distributed in an arbitrary way to the $r+s$ boundary circles. If $S$ is a \textit{closed} nonempty connected surface choose a Morse function on $S$ as above and consider a nonempty preimage of a regular value. Then the corresponding $R$-module map $R\rightarrow R$ assigned to $S$ is computed by breaking
the surface and the corresponding morphism $R\rightarrow A^{\otimes r}\rightarrow R$ where $r\geq 1$ is the number 
of circles at the regular value. The extended topological quantum field theory now inserts $\mu_a\otimes Id^{r-1}: A^{\otimes r}\rightarrow A^{\otimes r}$ between the two $R$-module maps. 
It is not hard to show that the above construction is well-defined.
A multiplication morphism can be moved through a connected surface using 
the fact that $\mu $ and $\Delta $ are $A$-bimodule maps.  
For example consider the pair of pants morphism $A\rightarrow A\otimes A$. If the pant is colored $a$ we associate to it the morphism 
$$\Delta \mu_a=(\mu_a\otimes Id)\Delta =(Id\otimes \mu_a)\Delta $$
Note that these identities follow from 
$\Delta (ba)=\Delta (b)a=\Delta (ba)=a\Delta (b)$
for all $b\in A$,
and express the commutativity of $\mu$ and the $A$-bimodule property of 
$\Delta $.
$\blacksquare$

\vskip .1in

Let $Cob_R(A)$ denote the category where the morphism sets of 
$Cob(A)$ are replaced by the free $R$-modules with bases the sets
of morphisms of $Cob(A)$. Then 
the above functor extends by $R$-linearity \textit{uniquely} to the functor:
$$F(\mathcal{F}): Cob_R(A)\rightarrow M(\mathcal{F}).$$
 
The definitions in the following section are motivated from analyzing the 
kernel of this functor $F(\mathcal{F})$, see section 5. 

\section{Skein modules of surfaces in $3$-manifolds.}  

Let $M$ be a $3$-manifold and let $\alpha \subset \partial M$ be 
a closed $1$-manifold. Let $\amalg $ denote disjoint union
and $|X|$ denote the number of components of a topological space $X$.  

\vskip .1in

Let $\mathcal{F}=(R,A,\Delta ,\varepsilon )$ be a Frobenius algebra as in section 2.
Choose a Frobenius system for $\mathcal{F}$ with $\Delta (1)=\sum_{i=1}^ru_i\otimes v_i$ and $u_i,v_i\in A$. 

\vskip .1in

For each subset $\mathfrak{a}\subset A$ let
$\mathcal{S}(M,\alpha, \mathfrak{a})$ denote the set of isotopy classes 
of properly embedded surfaces $S$ in $M$ with boundary $\alpha $ (isotopy relative to the boundary), and components colored by elements of $\mathfrak{a}$. 
The elements of $\mathcal{S}(M,\alpha , \mathfrak{a})$ are called \textit{$\mathfrak{a}$-colored surfaces} in $M$ bounding $\alpha $. For $\alpha =\emptyset $ we also consider the empty surface $\emptyset $ (with empty color) as an element of $\mathcal{S}(M, \emptyset ,\alpha)$. 

The most important cases are $\mathcal{S}(M,\alpha ):=\mathcal{S}(M,\alpha ,A)$, often just denoted $\mathcal{S}$, 
and $S(M,\alpha ,1)$. The second set will always be identified with the \textit{usual} 
set of isotopy classes of surfaces in $M$ bounding $\alpha $.
Formally, elements of $\mathcal{S}(M, \alpha)$ are pairs $(S,w)$, where $w\in H^0(S;A)$ is a locally constant continous map $w: F\rightarrow A$.

For each multplicatively closed $\mathfrak{a}\subset A$ we can also consider 
several colors on a component. Then the commutative
product of $A$ associates a unique element of $\mathcal{S}(M, \alpha ,\mathfrak{a})$.
 There is the map defined by forgetting colors:
$$\mathfrak{f}: S(M, \alpha ,A)\rightarrow S(M,\alpha ,1),$$
and
$\mathfrak{f}^{-1}(S)$ is in one-to-one correspondence with
$A^{|S|}$ for $S\in \mathcal{S}(M,\alpha ,1)$, at least up to possible 
reordering of 
components by isotopy. 

In order to describe relations we use the following local \textit{patch} notation. 
We let $(a)$ represent a local patch of a surface colored by $a\in A$.
This can be the color of a whole component, or there can be other colors on the same component, depending on the situation. Also in $(a)(a')$ the two patches may be on the same or distinct components. 
Sometimes we will also
allow additional colors outside of the surface patch.
Thus the component of the patch $(a)$ can be colored by a product of $a$ with 
other elements of $A$.  

Consider the $R$-submodule $\mathfrak{R}:=\mathfrak{R}(M,\alpha ;\mathcal{F})$ 
of $\mathcal{S}$ that is generated by the following three types of elements.

\noindent (1) (\textit{$R$-linearity}) For all $a_1,a_2\in A$ we have 
$$(a_1+a_2)-(a_1)-(a_2)\in \mathfrak{R},$$
and for $a\in A$ and $r\in R$ we have 
$$(ra)-r(a)$$ 
(Notice the similarity with the definition of tensor products).

\noindent (2) (\textit{sphere relations}) Let $(a)$ be a surface in $\mathcal{S}$ where the component colored by $a$ is a $2$-sphere bounding a $3$-ball in $M$. Then the difference between $a$ and the multiplication of the surface we get by omitting the sphere component by $\varepsilon (a)$ is in $\mathfrak{R}$.

\noindent (3) (\textit{neck cuttings}) Let $\gamma $ be a simple closed curve in the interior of a surface $S$ 
such that $\gamma $ bounds a disk $D$ in $M$ with $D\cap F=\gamma $. 
(Such a curve is automatically two-sided. This remark is due to 
Charlie Frohman.)
Then
$$(a)-\sum_{i=1}^r(au_i)(v_i)\in \mathfrak{R}.$$
Here the surface on the right hand side results by replacing
the annular neighborhood of $\gamma $ by two embedded disks in $M$,
which are the indicated pieces.
 
If $\gamma $ is separating the two disks are on different 
components and we say this is \textit{separating neck cutting}.
If $\gamma $ is non-separating then 
$$(a)-(\mu \Delta (a))$$
is in $\mathfrak{R}$ where the right hand surface is defined by cutting the 
handle, and we have \textit{nonseparating neck cutting}. 

\vskip .1in

We want to prove that the submodule $\mathfrak{R}$ defined as above 
only depends on the 
Frobenius algebra and not on the Frobenius system giving the elements
$\{(u_i,v_i)\}$ used in (3). In order to see this consider colored surfaces in $\mathcal{S}$
equipped with two disks embedded in the surfaces. Let $\widetilde{\mathcal{S}}$ denote the 
resulting set of isotopy classes. There is defined a mapping 
$$(A\otimes A)\times \widetilde{\mathcal{S}}\rightarrow C',$$
where $C'$ is the quotient of $R\mathcal{S}$ by the $R$-linearity 
relations,
as follows. Given $\sum_{j=1}^sx_i\otimes y_i\in A\otimes A$ and 
$S$ with the distinguished disks colored 
$a$ and $1$, define the image by
$$\sum_{i=1}^s(a_1x_i)(a_2y_i).$$
where we use the local notation as before.
Note that this is well-defined in $C'$ and only depends on 
the element $\sum_{i=1}^sx_i\otimes y_i\in A\otimes A$.
Now a curve $\gamma $ on a surface like in the neck-cutting relation 
defines several liftings to $\widetilde{\mathcal{S}}$
differing by the order of disks and 
the different ways of splitting $a=a_1a_2$ if $a$ is the coloring 
of the component containing $\gamma $. 
But because of the cocommutativity of $\Delta $ and the $A$-bimodule 
property the resulting element of $\mathcal{C}'$ does not depend on 
the choice of lifting.  
 
Thus $\mathfrak{R}$ is determined by the Frobenius algebra, and we have the following:

\vskip .1in

\noindent \textbf{Definition 4.1.} For each Frobenius algebra
$\mathcal{F}$ and $3$-manifold $M$ with closed $1$-manfold
$\alpha \subset \partial M$ there is defined the \textit{skein  module} 
$$C(M, \alpha ; \mathcal{F}):=R\mathcal{S}(M,\alpha ,A)/\mathfrak{R}(M,\alpha ;\mathcal{F}).$$
 Usually we omit
the \textit{coefficient system} from the notation. 
Let $C(M)=C(M;\mathcal{F})$ if $\alpha =\emptyset $.

\vskip .1in

The image of $(S,w )$ in $C(M,\alpha )$ is usually also denoted $(S,w)$, and 
we write $S$ for $(S,1)$, a surface with all components colored by $1\in A$.
We just write $1$ instead of $\emptyset$ in linear combinations. 

\vskip .1in

\noindent \textbf{Remarks 4.2.} (a) Obviously $\mathcal{S}(M,\alpha ,\mathfrak{a} )\neq \emptyset $ for $\mathfrak{a}\neq \emptyset $ if and only if $\mathcal{S}(M, \alpha ,1)\neq \emptyset$. 
If $\mathcal{S}(M,\alpha ,A)=\emptyset$ then $C(M,\alpha )=0$. But it is not a priori clear whether $\mathcal{S}(M,\alpha ,1)\neq \emptyset $ implies $C(M,\alpha  )\neq 0$. We will discuss some of these issues in section 7. 

\noindent (b) The reason to allow nonorientable surface in the definition of skein modules is motivated from Khovanov theory on surfaces \cite{M}, \cite{TT}.
Of course it is easy to define orientable or oriented versions of skein modules because
the relations (1)-(3) above preserve orientations. There are obvious definitions of those skein module, for each Frobenius algebra. For example let
$\tilde{C}(M,\ \tilde{\alpha} )$ be the quotient of the free $R$-module generated by
oriented properly embedded surfaces in $M$ bounding the oriented $1$-manifold 
$\tilde{\alpha }\subset \partial M$ by the above relations 
(1)-(3). 

\noindent (c) Given $S\in S(M, \alpha ,1)$ then the quotient of $R\mathfrak{f}^{-1}(S)$ by the submodule spanned by relations (1) is a certain quotient of $A^{\otimes |S|}$, the quotient coming from some action of the symmetric group corresponding to possible reordering of colored components by isotopy.
The structure of skein modules of surfaces in $3$-manifolds as above very 
much is determined by this 
algebraic similarity.

\noindent (d) Let $\mathcal{F}$ be graded such that $\varepsilon $ has degree
$-2$ and $\Delta $ has degree $2$. Then the skein module $C(M, \alpha )$ is graded using $deg(S,w)=-\chi (S)+\sum_jdeg(a_j)$, where the sum is over all components of the colored surface $S$ with $w$ assigning $a_j$ to the $j$-th component.

\noindent (e) Let $y\in R$ be invertible. Then multiplication by $y$ defines 
an isomorphism of modules 
$$C(M, \alpha ;\mathcal{F})\rightarrow C(M, \alpha ,\mathcal{F}^y).$$
But if $y\in A\setminus R$ is invertible the situation is more difficult 
because in general the skein modules are \textit{not} $A$-modules. 
For $y\in A$ invertible there is defined the natural bijective map:
$$\mathcal{S}(M, \alpha )\rightarrow \mathcal{S}(M, \alpha )$$
by $(S,w)\mapsto (S,aw)$, i.\ e.\ changing the colors of all components.
This is easily seen to be compatible with all relations except separating neck cuttings. 

\vskip .1in

Define $C_1(M,\alpha )\subset
C(M,\alpha )$ as the submodule generated by the image of
$S(M,\alpha ,1)$ in $C(M, \alpha )$. 
 
\vskip .1in

\noindent \textbf{Proposition 4.3.} \textit{If the powers of $\mu (\Delta (1))\in A$ form a set of $R$-generators of the algebra
$A$ then $C_1(M,\alpha )=C(M,\alpha )$.} 

\vskip .1in

\noindent \textit{Proof.} This follows because deleting a handle from a 
component of a surface in $M$ using a neck cutting amounts to multiplication of the color of that component by $\mu (\Delta (1))$. Thus we can eliminate the colors of components by adding trivial handles. $\blacksquare$ 

\vskip .1in

If 4.3 applies we say that the skein modules are \textit{geometric}. 
Usually 4.3. does not apply, even in the Bar-Natan case with $R=\mathbb{Z}$.  
$\mu (\Delta (1))=2X$. The powers of this element oviously do not generate 
$A$. If we replace $R$ by $\mathbb{Z}[\frac{1}{2}]$ and correspondingly $A$
the skein modules become geometric. For the universal rank two Frobenius
algebra $\mu (\Delta (1))=2X-h$. Again, if we replaced $R=\mathbb{Z}[h,t]$ by
$\mathbb{Z}[\frac{1}{2},h,t]$ then the module becomes geometric. But then
using some easy isomorphism one can see that the Frobenius 
algebra will just be the Gad Naot system (with an 
additional parameter $h$ added), see \cite{K1}, Examples 1. 

\vskip .1in

We will discuss the relation of the definition above with the examples 
previously discussed in the literature in section 8.

\vskip .1in

The following is immediate from the definition.

\vskip .1in

\noindent \textbf{Proposition 4.6.} Suppose that the $3$-manifold $M$ has $r$ components $M_i$, $i=1,\ldots ,r$ and
$(M,\alpha )=\amalg_{1\leq i\leq r}(M_i,\alpha _i)$,
$\alpha_i \subset \partial M_i$ for $i=1,\ldots ,r$. Then 
$$C(M,\alpha )=\bigotimes_{i=1}^r\mathcal{C}(M_i,\alpha_i).$$

\section{$2$-dimensional toplogical quantum field theory and \\
 \qquad skein modules.}

Consider the skein modules
$$\mathcal{C}(r,s):=C(D^2\times I,u^r\times \{1\}\amalg u^s\times \{0\}),$$
where $u^r$ respectively $u^s$ denotes the union of $r$ respectively $s$ 
standardly embedded trivial circles in $D^2$. 

There are $R$-module maps from
$\mathcal{S}(D^2\times I,ru\times 1\amalg su\times 0,A)$
onto the sets of morphisms of the category $Cob_R(A)$. 
The relations in the definition of skein modules can be 
applied \textit{abstractly} to the morphism sets of 
$Cob_R(A)$ to define a quotient category $Cob(\mathcal{F})$.
The objects of this category are still the objects of $Cob$
but the morphism sets are the quotients of the morphism modules 
of $Cob_R(A)$ by $R$-linearity and abstract sphere- and neck cutting 
relations. The resulting quotient category is denoted 
$Cob(\mathcal{F})$.

So there are induced $R$-module maps
$$\mathcal{C}(r,s)\rightarrow Cob(r,s)$$ 
where $Cob(r,s)$ is the morphism module from $r$ to $s$ circles
in the category $Cob(\mathcal{F})$. 
We want to prove that these are in fact isomorphisms.

\vskip .1in

\noindent \textbf{Theorem 5.1.} \textit{The functor $F(\mathcal{F})$ defined 
in} theorem 3.2., more precisely the extension to $Cob_R(A)$, 
\textit{factors through $Cob(\mathcal{F})$ and defines an isomorphism of categories
$$Cob(\mathcal{F})\rightarrow M(\mathcal{F}).$$}

\noindent \textit{Proof.} The observation that the functor factors through 
the quotient category $Cob(\mathcal{F})$ is by construction. For example 
a sphere colored $a\in A$ will map to the $R$-morphism $R\rightarrow R$ given by $\varepsilon \circ \mu_a \circ \iota =\varepsilon (a)Id_R$. But this is also the image of 
$\varepsilon (a) \emptyset $, where $\emptyset $ is the empty surface
(which maps to the identity $Id_R: R\rightarrow R$).
A similar argument using equation (2) from section 2 proves that 
the neck cutting relations are in the kernel of the topological 
quantum field theory. Thus 
$$F(\mathcal{F}): Cob(\mathcal{F})\rightarrow M(\mathcal{F})$$
is well-defined and onto by the definition of $M(\mathcal{F})$. 
We have to prove injectivity on the set of morphisms.
For $r,s\geq 0$ 
consider the $R$-morphism 
$$\phi_{r,s}: A^{\otimes (r+s)}\rightarrow Cob(r,s)$$ 
defined by mapping $a_1\otimes a_r\otimes b_1\otimes b_s$ to the 
surface consisting of $r+s$ disks colored with $a_1,\ldots ,b_s$.
This is well-defined in $Cob(r,s)$ because of the 
$R$-linearity relations. Consider the composition with
$F(\mathcal{F})$. This is the $R$-morphism: 
$$A^{\otimes (r+s)}\rightarrow Hom(A^r,A^s)$$
$$a_1\otimes \ldots \otimes a_r\otimes b_1\otimes \ldots \otimes b_s \mapsto 
(a'_1\otimes \ldots \otimes a'_r\mapsto \varepsilon (a_1a'_1)\ldots \varepsilon (a_ra'_r)b_1\otimes \ldots \otimes b_s).$$
It is not hard to see that these $R$-morphisms are injective. 
To keep notation simple we show the case $r=2, s=1$. 
Let $\sum_{k=1}^Na_k\otimes b_k\otimes c_k$ 
be in the kernel of the composition. It is 
mapped to the $R$-morphism 
$A\otimes A\rightarrow A$ defined by 
$$a \otimes b \mapsto 
\sum_{k=1}^N\varepsilon (aa_k)\varepsilon (bb_k)c_k,\leqno(4)$$
which we assume to be the trivial map. 
Now for $k=1,\ldots ,N$ we can write using equation (2) of section 2:
$$a_k=\sum_{i_k=1}^r\varepsilon (a_ku_{i_k})v_{i_k}$$
and
$$b_k=\sum_{j_k=1}^r\varepsilon (b_ku_{j_k})v_{j_k}.$$
Thus 
$$\sum_{k=1}^Na_k\otimes b_k\otimes c_k=\sum_{k=1}^N\sum_{i_k=1}^r\sum_{j_k=1}^r\varepsilon (a_ku_{i_k})\varepsilon (b_ku_{j_k})v_{i_k}\otimes v_{j_k}\otimes c_k.$$
But substituting $a=u_{i_k}$ and $b=u_{j_k}$ for fixed values of $i_k,j_k\in \{1,\ldots ,r\}$ we get
$$\sum_{k=1}^N\varepsilon (a_ku_{i_k})\varepsilon (b_ku_{j_k})c_k=0$$
and thus the sum is vanishing. $\blacksquare$

\vskip .1in

It follows from the proof of 5.1. that the $R$-module maps
$\phi_{r,s}$ are isomorphisms. This proves the second half of:

\vskip .1in

\noindent \textbf{Corollary 5.2.} \textit{For all $r,s$ we have
natural isomorphisms: 
$$\mathcal{C}(r,s)\cong Cob(r,s)\cong A^{\otimes (r+s)}.$$
Under these isomorphisms the composition morphisms
$$Cob(s,t)\otimes Cob(r,s)\rightarrow 
Cob(r,t)$$
correspond to 
$$(b'_1\otimes b'_s\otimes c_1\otimes c_t)\otimes (a_1\otimes a_r\otimes b_1\otimes b_s)\mapsto \varepsilon (b_1b'_1)\ldots \varepsilon (b_sb'_s)(a_1\otimes a_r\otimes c_1\otimes c_t).$$
}

\noindent \textit{Proof.} Since surfaces embedded in $D^2\times I$ are 
completely compressible we can apply neck cutting and sphere 
relations to see that the modules
$\mathcal{C}(r,s)$ also are also $R$-generated by 
$A$-colored disks bounding the corresponding input and output circles.
So the $R$-module maps:
$$\psi_{r,s}: A^{\otimes (r+s)}\rightarrow \mathcal{C}(r,s)$$
defined by coloring disks as above is onto, and 
factors the $R$-module map 
$\phi_{r,s}$ from above. Because $\phi_{r,s}$ is injective, 
$\psi_{r,s}$ is also injective. Thus the projections 
$\mathcal{C}(r,s)\rightarrow Cob(r,s)$ are isomorphisms. $\blacksquare$ 

\vskip .1in

Theorem 5.1. and its corollaries contain several results concerning 
Bar-Natan modules previously discussed in the literature, see e.\ g.\
\cite{N1}, \cite{N2}. For example the case $r=s=0$ gives the following result:

\vskip .1in

\noindent \textbf{Corollary 5.3.} \textit{For each Frobenius algebra 
$\mathcal{F}$
the skein module $C(D^3,u^r)$ is isomorphic to $A^{\otimes r}$.
In paricular for $r=0$ we have 
$$C(D^3)\cong R.$$
The image of a connected surface of genus $g$ colored by $a\in A$ is
$$\varepsilon ((\mu (\Delta (1))^ga)$$
and the isomorphisms maps disjoint unions to products.} $\blacksquare$ 

\section{Some properties of skein modules.}

The following result has been proved for Nar-Natan modules in \cite{FK}. 
For $i=1,2$ let $M_i$ be two $3$-manifolds with 
disks $D_i\subset \partial M_i$ and
closed $1$-manifolds 
$\alpha_i\subset \partial M_i\setminus D_i$.
Then we can define the boundary connected sum $M:=M_1\sharp M_2$ with the 
closed $1$-manifold $\alpha :=\alpha_1\cup \alpha _2$ in its boundary. 
Note that $M_i\subset M_1\sharp M_2$ for $i=1,2$.

\vskip .1in

\noindent \textbf{Theorem 6.1.} 
$$C(M_1\sharp M_2,\alpha_1\cup \alpha_2)\cong C(M_1,\alpha_1)\otimes C(M_2,\alpha_2).$$ 

\noindent \textit{Proof.} 
Let $D\subset M_1\sharp M_2$ be the disk along the connected sum has
been formed.
First define a $R$-module map:
$$\phi: C(M_1,\alpha_1)\otimes C(M_2,\alpha_2)\rightarrow C(M,\alpha ).$$
Since $C(M_1,\alpha_1)\otimes C(M_2,\alpha_2)$ is $R$-generated by 
the image of $S(M_1,\alpha_1,A)\otimes S(M_2,\alpha_2,A)$ under 
the tensor product of the natural projections,
we can let $\phi $ map the image of $S_1\otimes S_2$ to the disjoint union
$S_1\amalg S_2\subset M_1\sharp M_2$, which bounds $\alpha $. It is immediate from the definitions 
that this map is well-defined.
This map is surjective because we can apply neck cutting relation to
all the circles in the transversal intersection $D\cap S$ for a surface 
$S\in \mathcal{S}(M, \alpha)$. 
Note that the set of intersection circles is partially ordered by 
nesting. The neck cuttings have to be performed starting from 
innermost circles. 
In order to show that $\phi $ is injective we will construct the 
inverse $R$-module map:
$$\psi : C(M,\alpha )\rightarrow C(M_1,\alpha_1)\otimes C(M_2,\alpha_2)$$
as follows. 
For each generating (colored) surface $S\subset M$ bounding $\alpha $ 
apply neck cuttings to the transverse intersection $S\cap D$.  

We have to prove that the resulting linear combination in
$C(M_1,\alpha_1)\otimes C(M_2,\alpha_2)$ does not depend on the 
choice of surface in its isotopy class relative boundary. 
Now isotopies $S_t$ where the intersection $D\cap S_t$ is 
only changes by isotopy in $D$ 
will give rise to isotopies of the cut surfaces in $M_1$ respectively 
$M_2$. Otherwise it suffices to
consider elementary isotopies with a single time parameter $t_0$ for 
which the surface $S_{t_0}$ is not transversal as follows (or its inverse): 
(i) ths surface touches the disk at some 
point, and thus pushing through the disk will create a new intersection 
circle, (ii) the surface intersects $D$ in a 
\textit{saddle}, so for $t<t_0$ the intersection consists of two circles 
that are connected to a single circle when pushing the saddle through $D$
from $M_1$ into $M_2$.
Because the cutting is applied following the partial order of circles we
can assume that all these circles are innermost.  
Now in case (i) there is a $2$-sphere bounding a $3$-ball in $M$
and application of the neck cutting does not contribute. 
In the second case we compare the result of the two neck cuttings for 
$t<t_0$ with the neck cutting along the one circle for $t>t_0$. 
But an additonal neck cutting in $M_2$ gives the same result.
What we have to check can, because of $R$-multilinearity, be expressed 
in the patch notation by:
$$\sum_{i,j}(au_iu_j)(v_j)(v_i)=\sum_{ij}(au_i)(u_j)(v_iv_j)$$
for all $a\in A$. 
But this is the relation
$$(Id\otimes \Delta )\circ \Delta =(\Delta \otimes Id)\circ \Delta ,$$
which follows from cocommutativity of $\Delta $ and $A$-bilinearity.
It corresponds to switching two consecutive saddle points of a surface 
in the topological quantum field theory associated to $\mathcal{F}$.   
Finally compatibility with relations (1)-(3) is easly seen from locality of the
relations.(For the neck cuttings one has to use the disks in $M$ bounding 
the cutting curves on $M$ to find isotopies to push the disks away from $D$.) 
 Thus $\psi $ is well-defined, and $\psi \circ \phi=Id$ proves that 
$\phi $ is also injective. $\blacksquare$  

\vskip .1in

The special case of $M_2$ a $3$-ball, with a trivial loop $u$ in its boundary
disjoint from the connected sum disk, proves by induction on $r$ the following:

\vskip .1in

\noindent \textbf{Corollary 6.3.} \textit{Let $r\geq 0$. Suppose $u^r$ is any disjoint union of  simple closed loops in $\partial M$, each of which bounds a disk in $\partial M$ and is contained in the complement of a $1$-manifold $\alpha \subset \partial M$. Then there is an isomorphism of $R$-modules
$$C(M,\alpha \amalg u^r)\cong C(M,\alpha )\otimes A^{\otimes r}.\qquad \qquad \blacksquare$$}

\noindent \textbf{Remark 6.4.} For $\mathcal{F}$ free one can also give 
give a direct proof of 6.2. without reference to 5.3.
(which in fact reproves 5.3.). Again it suffices to prove the case $r=1$.
Let $\mathfrak{b}=\{u_i|i=1,\ldots ,r\}$ be a $R$-basis of $A$ and 
$\{w_j|j=1,\ldots ,r\}$ be a dual basis such that 
$\varepsilon (u_iw_j)=\delta_{ij}$.
Using 8.1. we can identify $C(M,\alpha )$ and $C(M,\alpha ,{\mathfrak{b}})$
(see 8.1.).
Then define the $R$-module map
$$\phi: \mathcal{S}(M,\alpha ,\mathfrak{b})\otimes A\rightarrow \mathcal{S}(M, \alpha \amalg u,\mathfrak{b})$$ 
by 
$$\phi((S,w)\otimes b)=  
(S',w')$$
for $b\in \mathfrak{b}$. 
Here $S'$ is the disjoint union of $S$ with a disk $D$ bounding $u$ 
in a collar neighborhood of $\partial M$ and define $w'|S=w$, 
$w'|D=b$. 
It is important to observe that
elements of $R$ can move freely between components of $S$ and the disk
using the $R$-linearity relations. 
Because the skein relations are local there is the induced $R$-module map 
$C(M,\alpha )\otimes A \rightarrow \mathcal{C}(M,\alpha \amalg u)$. This map is onto because we can apply the neck cutting relations to elements of $\mathcal{S}(M,\alpha \amalg u,\mathfrak{b})$ to expand in linear combinations of $\mathfrak{b}$-colored surfaces given by a marked disk 
bounding $u$ and a marked surface bounding $\alpha $.
It suffices to prove that $\phi $ is injective. 
Let $\phi (\sum_{i=1}^rx_i\otimes u_i)=0$  
with $x_i\in C(M, \alpha )$ for 
and $i=1,\ldots ,r$.
Then define for $j=1,\ldots ,r$:
$$\psi_j: C(M,\alpha \amalg u )\rightarrow C(M,\alpha )$$
on marked surfaces $(S,w)\in C(M, \alpha \amalg u)$ 
by adding a collar to $M$ along $\partial M$, capping off $u$ by
a disk colored $w_j$ and adding a cylinder over $\alpha $. This 
results in a well-defined element of $C(M,\alpha )$.  
It follows from the definitions that for $j=1,\ldots ,r$
$$\psi_j(\phi (\sum_{i=1}^rx_i\otimes u_i))=\sum_{i=1}^rx_i\varepsilon (u_iw_j)=x_j$$
Thus $\phi $ is also injective.

\vskip .1in

The proof of the following result is very similar to the proof of 
6.3.

\vskip .1in

\noindent \textbf{Proposition 6.5.} \textit{For $(M, \alpha )$ let 
$(M', \alpha )$ be the result of attaching a $1$-handle $H$ to $M$ along 
$\partial M$ such 
that 
$H\subset \partial M\setminus \alpha$. Then the inclusion 
$M\subset M'$ induces the isomorphism
$$C(M,\alpha )\rightarrow C(M',\alpha ).\quad \blacksquare$$}

\noindent \textbf{Corollary 6.6.} 
\textit{$C(H_g)\cong R$ for $H_g$ the genus $g$ handlebody.} $\blacksquare$

\vskip .1in

It follows that $C(H_g,u^r)\cong A^r$ for trivial curves in $\partial H_g$.   
Also $C(H_g,\mu )\cong A$ for a meridian $\mu $ follows by attaching a 
$1$-handles 
to $(D^3,u)$ suitably.

\vskip .1in

In general it seems difficult but interesting to determine 
$C(\Sigma \times I, \alpha )$ for $\alpha \subset \partial (\Sigma \times I)$, see also 7.5.

\section{Naturality properties of the skein modules.}

First we discuss a universal coefficient theorem. 
Let $\mathcal{F}=(R,A,\Delta,\varepsilon)$ and $\mathcal{F'}=(R',A',\Delta',\varepsilon')$ be two Frobenius algebra. Let $\phi : \mathcal{F}\rightarrow \mathcal{F}'$ be a \textit{morphism of Frobenius algebras}. This an $R$-algebra morphism
$A\rightarrow A'$, also denote $\phi $, which commutes with all obvious diagrams relating the structure maps of $\mathcal{F}$ with those of $\mathcal{F}'$. For example $\iota '\circ \phi |R=\phi \circ \iota$, or
$$(\phi \otimes \phi )\circ \Delta = \Delta '\circ \phi$$ 
Note that $\phi $ induces on $R'$ the structure of $R$-module.
A morphism of Frobenius algebras is called \textit{epimorphism} if both morphisms 
$A\rightarrow A'$ and $R\rightarrow R'$ are onto. An epimorphism of Frobenius algebras is called \textit{coefficient induced} if 
$A'=A\otimes_R R'$ and all structure maps of $(R',A',\Delta ',\varepsilon ')$ are induced via tensor products $\otimes_R$ in the obvious way from $(R,A,\Delta ,\epsilon )$.

\vskip .1in

\noindent \textbf{Example 7.1.} Let $G,G'$ be two finite groups. 
Then each group epimorphism $\phi : G\rightarrow G'$ extends to a 
epiphism of Frobenius algebras 
$\mathcal{F}_G\rightarrow \mathcal{F}_{G'}$
defined by the natural morphism of group algebras.
In particular this applies to $G\rightarrow \{1\}$.

\vskip .1in 

The first claim of following proposition is obvious, and the second one
is proven in the very same way as the universal coefficient
theorem for link skein modules, see \cite{P}.

\vskip .1in

\noindent \textbf{Proposition 7.2.} Let $\phi : \mathcal{F}\rightarrow 
\mathcal{F'}$ be an epimorphism of Frobenius algebras. Then $\phi $ induces an epimorphism of $R$-modules
$$C(M,\alpha ,\mathcal{F})\rightarrow C(M,\alpha ,\mathcal{F}'),$$
where $C(M,\alpha ,\mathcal{F}')$ is considered as $R$-module  
via $R\rightarrow R'$.
If $\phi $ is a coefficient induced epimorphism then there is a natural 
isomorphism
$$C(M,\alpha ;\mathcal{F'})\cong C(M,\alpha ;\mathcal{F})\otimes R'.\qquad \qquad \blacksquare$$  

\vskip .1in

Consider for $\alpha \subset \partial M$ the exact sequence:
$$
\begin{CD}
H_2(\alpha ;\mathbb{Z}_2)@>>>H_2(M;\mathbb{Z}_2)@>>>H_2(M,\alpha ;\mathbb{Z}_2)@>\partial >>H_1(\alpha ;\mathbb{Z}_2 )@>j_*>>H_1(M;\mathbb{Z}_2)
\end{CD}
$$
and 
$$\partial ^{-1}[\alpha ]\subset H_2(M, \alpha ;\mathbb{Z}_2)=H_2(M, N(\alpha );\mathbb{Z}_2),$$
where $[\alpha ]$ is the \textit{fundamental class} and $N(\alpha )$ is a neighborhood of $\alpha $ in $\partial M$. 

By the Pontrjagin-Thom construction the set 
$$\partial^{-1}[\alpha ]\subset H_2(M, N(\alpha );\mathbb{Z}_2)\cong H^1(M/N(\alpha );\mathbb{Z}_2)\cong [M/N(\alpha ),\mathbb{R}P^{\infty }]$$ 
can be interpreted as the set of bordism classes of surfaces in $M$ with boundary $\mathbb{Z}_2$-\textit{homologous} to $\alpha $. It is possible that $\partial^{-1}[\alpha ]=\emptyset $, in 
which case $S(M, \alpha , A)=\emptyset $ and $C(M, \alpha )=0$. (For example 
take for $\alpha $ an odd number of longitudes in the boundary of a solid torus.)
On the other hand, each element in $\mathcal{S}(M, \alpha ,A)$ determines 
an element in $H_2(M, N(\alpha );\mathbb{Z}_2)$.  

\vskip .1in

The skein modules of surfaces always decompose with respect to 
$\mathbb{Z}_2$-homology.

\vskip .1in

\noindent \textbf{Proposition 7.3.}
$$C(M,\alpha )\cong \bigoplus_{a\in \partial^{-1}[\alpha ]}C_a(M,\alpha ),$$
where $C_a(M,\alpha )$ is the skein module spanned by all surfaces 
representing $a\in H_2(M, \alpha :\mathbb{Z}_2)$. 
If b$j_*[\alpha]=0$ then $\partial^{-1}[\alpha ]$ is in one-to-one correspondence with $H_2(M;\mathbb{Z}_2)$. Otherwise $\partial^{-1}[\alpha ]=\emptyset$.

\vskip .1in

\noindent \textit{Proof.} This follows from the fact that the skein relations
imply $\mathbb{Z}_2$-homology of the underlying surfaces and the exactness of the above homology sequence. $\blacksquare$   

\vskip .1in

For skein modules of surfaces bounding oriented curve systems the analogous result holds with $\mathbb{Z}$-homology replacing $\mathbb{Z}_2$-homology. The discussion of skein modules of orientable surfaces bounding curve systems requires to study the image of $\mathbb{Z}$-homology in $\mathbb{Z}_2$-homology.

\vskip .1in

Suppose $\alpha $ is two-sided on $\partial M$ such that $N(\alpha )$ is a union of tori. Then any surface with boundary $\mathbb{Z}_2$-homologous to $\alpha $ can be modified by glueing annuli to pairs of parallel boundary components 
to a surface bounding $\alpha $. In for each $a\in \partial^{-1}[\alpha ]$ there exists an element of $S(M, \alpha ,A)$ with homology class $\alpha $.

\vskip .1in

\noindent \textbf{Proposition 7.4.} \textit{Suppose $S_1, S_2\in S(M, \alpha )$ are two surfaces 
with $[S_1]=[S_2]\in H_2(M,\alpha ;\mathbb{Z}_2)$. Then by adding tubes to $S_1$ and 
$S_2$ we get surfaces $S_1'$ and $S_2'$ which are isotopic relative boundary.}

\vskip .1in

\noindent \textit{Proof.} It obviously suffices to consider the case $M$ connected. The basic idea is classical and often used in the discussion of $S$-equivalence of Seifert surfaces (see \cite{Ka} for a detailed account) with the main argument given in \cite{KL}. Actually a quite subtle argument is necessary to prove first that $S_1$ and $S_2$ bound a $3$-manifold 
$W\subset M\times I$ such that $W\cap (\partial M\times I)=\alpha \times I$.  
This is shown by obstruction theory of maps into $\mathbb{R}P^{\infty}$.
Then by adding $1$-handles, i.\ e.\ tubes, we can make both $S_1, S_2$ connected. A handle cancellation argument, see \cite{KL}, then proves the claim.
$\blacksquare$

Let $d: M\rightarrow M$ be a diffeomorphism of $M$. Then $d$ induces a $R$-module map:
$$d_*: C(M, \alpha )\rightarrow C(M, h(\alpha ))$$
The induced map is not changing under diffeomorphisms that are constant on the boundary. In particular the modules $C(M, \alpha )$ are representations of the \textit{mapping class group} $Diff(M, \partial M)$ of isotopy classes of diffeomorphisms that are the identity map on $\partial M$ 

\vskip .1in

\noindent \textbf{Example 7.5.} For the understanding of Khovanov theory of diagrams on oriented surfaces $\Sigma $ the calculation of the skein modules
$$C(\Sigma \times I, \alpha \times \{1\} \cup \beta \times \{0\})$$
for $\alpha ,\beta \subset \Sigma $ two given $1$-dimensional closed manifolds is important, at least for the universal rank $2$ Frobenius algebra. By 6.3. we can assume that both $\alpha $ and $\beta $ are essential. The skein modules here are also morphism sets of a certain category $Cob(\Sigma ,\mathcal{F})$, essentially first defined by Turner and Turaev in \cite{TT} 
(and denoted $\mathcal{U}Cob(\Sigma )$ there to emphasize the possibility of unorientable surfaces).
in some case. The understanding of the $R$-algebras
$$C(\Sigma \times I, \alpha \times \{0,1\})$$ 
is interesting for the understanding of the \textit{tautological} topological quantum field theories as suggested by Bar-Natan. 
The important cases here are the calculations of the two algebras 
$\alpha $ a simple closed non-separating curve 
and for $\alpha $ empty. See \cite{FK} for the calculation in the Bar-Natan 
case.

\section{Presentations of skein modules from presentations of $A$.}

First we show how a presentation of the $R$-module $A$ induces a presentation 
of the $R$-module $C(M, \alpha )$.

\vskip .1in

Let $\mathfrak{g}\subset A$ be a generating set of $A$ and 
$\mathfrak{r}\subset R\mathfrak{g}$ be a set of relations such that
$$A\cong R\mathfrak{g}/span(\mathfrak{r}),$$
where $span(\mathfrak{r})$ is the submodule of $R\mathfrak{g}$ generated by 
$\mathfrak{r}$. Let $\mathfrak{p}$ denote the \textit{presentation} 
$(\mathfrak{g}, \mathfrak{r})$. 

\vskip .1in

Recall that
$$\mathcal{S}(M, \alpha ,\mathfrak{g})\subset \mathcal{S}(M,\alpha , A)$$
is the set of isotopy classes of surfaces in $M$ bounding $\alpha $ with colors of all components in $\mathfrak{g}$. If a component is colored 
by several elements of $A$ then we always assume that their product is in 
$\mathfrak{g}$, even if $\mathfrak{g}$ is not multiplicatively closed.
The elements of $S(M, \alpha ,\mathfrak{g})$ are called
\textit{marked surfaces} with respect to $\mathfrak{g}$. 
Then let 
$$\mathfrak{R}(M,\alpha ,\mathfrak{p}) \subset R\mathcal{S}(M, \alpha ,\mathfrak{g})$$
to be the submodule generated by (using the patch notation from section 4):

\noindent (1') $\sum_{i=1}^s (g_i)=0$ where $\sum_{i=1}^s g_i\in \mathfrak{r}$
and $g_i\in \mathfrak{g}$ for $i=1,\ldots ,s$.

\noindent (2') sphere relations for all $2$-spheres bouding $3$-ball which are colored by elements of $\mathfrak{g}$.

\noindent (3') neck cutting relations where we assume that the component of the neck is colored by some element of $A$ and the right hand side is expanded 
using $R$-linearity into a linear combination of $\mathfrak{g}$-colored surfcaes.

\vskip .1in
 
Let $C(M,\alpha ;\mathfrak{p})$ denote the quotient of $R\mathcal{S}(M,\alpha ,\mathfrak{g})$ by $\mathcal{R}(M,\alpha ,\mathfrak{p})$. 

\vskip .1in

\noindent \textbf{Theorem 8.1.} \textit{The inclusion 
$$\mathcal{S}(M,\alpha ,\mathfrak{g})\subset \mathcal{S}(M,\alpha ,A)$$
induces the isomorphism
$$C(M,\alpha ,\mathfrak{p})\cong C(M,\alpha ).$$}

\noindent \textit{Proof.} Obviously $\mathfrak{R}(M, \alpha ,\mathfrak{p})$ 
is contained in the kernel of the induced $R$-module map
$$R\mathcal{S}(M,\alpha ,\mathfrak{g})\rightarrow C(M,\alpha )$$
which is onto because $\mathfrak{g}$ is a generating set and the 
$R$-linearity relations hold in $C(M, \alpha )$. 
Let $\phi $ be the $R$-module map induced on the quotient.
Construct an inverse of $\phi $ by first defining 
$$\mathcal{S}(M,\alpha ,A)\rightarrow C(M,\alpha ,\mathfrak{p})$$
using expansion of elements of $A$ in terms of $\mathfrak{g}$ and 
$R$-linearity. It is clear from the definitions that 
$\mathfrak{R} =\mathfrak{R}(M,\alpha ;\mathcal{F}) \subset R\mathcal{S}(M,\alpha ,A)$ is contained in the kernel of 
this mapping. The induced $R$-module map shows that $\phi $ is an isomorphism.
$\blacksquare$

\vskip .1in

\noindent \textbf{Remarks 8.2.} (a) Suppose $\mathfrak{g}$ is closed under 
multiplication in $A$. Then the neck cutting relations $(3')$ above are just 
special neck cutting relations (3). 

\noindent (b) If $\mathfrak{g}$ is a basis and $\mathfrak{r}=0$ then there are no relations $(1')$.

\vskip .1in

\noindent \textbf{Examples 8.3.} (a) Consider the 
\textit{universal rank $2$ Frobenius algebra} $\mathcal{F}_{\mathbf{U}}$ with
$\mathfrak{g}=\{X^n|n=0,1,\ldots \}$. Then $\mathcal{S}(M,\alpha,\mathfrak{g})$ is the set of dotted surfaces where a dot on a component corresponds to a 
color $X$. The relations $(1')-(3')$ above in this case have been used
in \cite{C} to define the \textit{universal Bar-Natan module}. 
But note that $A$ is free with basis $\{1,X\}$. 

\noindent (b) If we reduce the universal rank $2$ Frobenius algebra by $t=h=0$ we have $R=\mathbb{Z}$ and the relation that each surface with a component with two dots is $0$ because $X^2=0$ (two dot relation).  
This is the Frobenius algebra $\mathcal{F}_{BN}$ studied in \cite{AF} and first described in 
\cite{BN}. It is relevant to the categorification of the Jones polynomial as originally described by Khovanov.
In this case the modules $C(M,\alpha )$ is the \textit{Bar-Natan module} of the $3$-manifold 
$M$ with respect to $\alpha  \subset \partial M$.
They have been computed for $\alpha =\emptyset $ in \cite{AF} for Seifert fibred $3$-manifolds.
   
\noindent (c) In \cite{N1}, \cite{N2} Gad Naot defines the following version of Bar-Natan modules. He takes the quotient of the free abelian group generated 
by the set of elements $\mathcal{S}(M,\alpha ,A)$ with all components 
colored $1$ (the usual isotopy classes of surfaces) by the Bar-Natan relations:
(2'') if a surfaces has a $2$-sphere component then the surface is trivial,
(3'') $2S=S_++S_-$, where $S$ is a surface with a neck given by a compressible
loop $\gamma $ on $S$, $S_{\pm }$ are the two surfaces which result by 
cutting as before but with adding an additional trivial $1$-handle inserted 
into the \textit{left} respectively \textit{right} hand disk.  
Then he proves that adding two $1$-handles into a component does not depend on the component. Even though his argument is given in the abstract setting it is easy to see that it can be generalized to arbitrary $3$-manifolds using connectedness and the possibility of unorientable surfaces. 
In this way the Gad-Naot module becomes a $\mathbb{Z}[T]$-module where $T$ is the two handles operator. Gad Naot discusses the situation if $\mathfrak{1}{2}$ is added to the ground ring. Then it is easy to see that the resulting module 
is essentially the $h=0$ reduction of our skein module for the universal
rank $2$ Frobenius algebra. We call the resulting Frobenius algebra
$\mathcal{F}_{GN}$.  
Note that labeling a component by $X$ corresponds to
multiplication by $\frac{1}{2}$ and adding a trivial $1$-handle to this component. Multiplication by $T$ in Gad Naots module thus corresponds to multiplication by $4t$ in our reduction module. 

\vskip .1in

The skein modules from \cite{BN}, \cite{N1} and \cite{N2} actually are originally based on the $4$-tube relation instead of the neck cutting relation. 
It has been proved in \cite{BN} that the neck cutting relation implies the $4$-tube relation. If $2$ is invertible in $R$, like in $\mathbb{Z}[\frac{1}{2},t]$ then the $4$-tube relation also imples the neck cutting relation. 
The $4$-tube relation is motivated from the Khovanov construction, giving the essential invariance of link homology under Reidemeister moves.
But from the viewpoint of topological quantum field theory or Frobenius algebras the neck cutting relations seems is more fundamental. 

\vskip .1in

Following 2.2. we calculate for the Gad Naot modules 
$\varepsilon_n=0$ for $n$ even nonnegative integers and 
$\varepsilon_n=t^{\frac{n-1}{2}}$ for 
$n$ odd integers. Also $\mu (\Delta (1))=2X$. 
Thus for a genus $g$ trivially colored surface 
by 5.3., 
the image in $C(D^3)$ is  
$$(1+(-1)^{g+1})^gt^{\frac{g-1}{2}}\in \mathbb{Z}[t],$$
see \cite{N1}.

Note that examples (a)-(c) define coefficient induced epimorphisms of Frobenius algebras
$$\mathcal{F}_{\mathbf{U}}\rightarrow \mathcal{F}_{GN}\rightarrow \mathcal{F}_{BN}.$$

\noindent (d) Let $\mathcal{F}_{\{e\}}$ be the Frobenius algebra defined in 
2.2 (b) for $G=\{e\}$ the trivial group. Then $A=R$ with identity $1=e$ and the neck cutting relation following $\Delta (1)=e\otimes e$ allows to replace a surface by the surface resulting from cutting the neck. Moreover, the trivially colored $2$-sphere is $1$ in this case. Thus it follws that
$$C(M, \alpha ;\mathcal{F}_{\{e\}})\cong \mathbb{Z}H_2(M;\mathbb{Z}_2)$$
for each $\alpha $ with $j_*[\alpha ]=0$.
If $j_*[\alpha ]\neq 0$ then $C(M, \alpha ;\mathcal{F}_{\{e\}})=0$, then  
using 7.1 and 7.2 it follows that $C(M, \alpha ;\mathcal{F}_G)\rightarrow \mathbb{Z}H_2(M;\mathbb{Z}_2)$ is an epimorphism.
Thus skein modules for group algebra Frobenius algebras are natural \textit{deformations} of the group algebra $\mathbb{Z}H_2(M;\mathbb{Z}_2)$.

\section{Presentations of skein modules from incompressible surfaces.}

First we generalize an important observation of Asaeda and Frohman \cite{AF}.
Recall that a \textit{compression disk} $D$ for a surface $S\subset M$ is a disk in $M$
such $D\cap S=\gamma $ where $\gamma $ is an essential curve on $S$.
Then $\gamma $ is called a \textit{compression curve}.  
The surface $S$ is called incompressible if there are no compression disks
and no components which are $2$-spheres bounding $3$-balls in $M$. 
A surface which does not have any compression curves but possibly has 
$2$-sphere components bounding $3$-balls in $M$ is called \textit{weakly incompressible}. Colored surfaces are called (weakly) incompressible if their 
underlying topological surfaces are (weakly) incompressible. For each $\mathfrak{a}\subset A$ let $\mathcal{I}(M, \alpha ,\mathfrak{a})$ denote the set of incompressible surfaces in $M$ bounding $\alpha $ with components colored by elements of $\mathfrak{a})$.
  
\vskip .1in

\noindent \textbf{Theorem 9.1.} \textit{Let $A$ be generated by $\mathfrak{g}\subset A$. Then the skein module $C(M,\alpha )$ is generated by the images
in $C(M, \alpha )$ of the set 
$\mathcal{I}(M,\alpha ,\mathfrak{g})\subset \mathcal{S}(M,\alpha ,\mathfrak{g})$ of marked incompressible surfaces.}

\vskip .1in

\noindent \textit{Proof.} By 8.1. the module is generated by marked surfaces. If a surface $F$ is not weakly incompressible then there is a compression curve $\gamma $ on $F$.
If $\gamma $ is not $2$-sided on $F$ then the normal bundle of $\gamma $ in $M$ is nontrivial. 
But this is not possible since $\gamma $ bounds a disk in $M$. (The fact
that orientability of $M$ is not relevant in this situation and the argument above, have been pointed out to me by 
Charlie Frohman.)
Thus we can apply a neck cutting relation and expand in terms of marked 
surfaces.  
After finitely many steps the surface will be weakly incompressible. 
Finally we can apply sphere relations to eliminate all $2$-spheres which bound $3$-balls in $M$.
$\blacksquare$  
 
\vskip .1in

Recall that a $3$-manifold is irreducible if each $2$-sphere in $M$ bounds a $3$-ball. It is a simple but important observation that for irreducible $M$ we can assume that neck cutting relations only apply to curves $\gamma $ which are essential on a colored surface $S$. In fact a neck cutting relation 
for an inessential loop on a surface is in patch notation:
$$(a)=\sum_i(au_i)(v_i)$$
with the second patch indicating a $2$-sphere, which bounds a $3$-ball in $M$
because of irreducibility.
But also
$$\sum_i(au_i)(v_i)
=\sum_i \varepsilon(v_i)(au_i)=(\sum_i\varepsilon(v_i)au_i)
=((Id\otimes \varepsilon )\circ \Delta (a))=(a)$$
using the $R$-linearity relations, sphere relations and equation $(1)$. 

\vskip .1in

Let $\mathfrak{p}=(\mathfrak{g}, \mathfrak{r})$  
be as in section 8. Let  
Let $\mathfrak{C}(M, \alpha ,\mathfrak{p})\subset  \mathfrak{R}(M, \alpha, \mathfrak{p})\cap R\mathcal{I}(M,\alpha ,\mathfrak{g})$ be the submodule of $R\mathcal{I}(M,\alpha ,\mathfrak{g})$, which is generated by all relations (1') with all colored surfaces in
$\mathcal{I}(M,\alpha ,\mathfrak{g})$.

We will next define the \textit{tubing submodule} 
$$\mathfrak{T}(M, \alpha ,\mathfrak{g}) \subset 
R\mathcal{S}(M, \alpha ,\mathfrak{g})$$
of tubing relations in $M$. The importance of this module comes from the following:

\vskip .1in

\noindent \textbf{Theorem 9.2.} \textit{Let $\mathcal{F}$ be a Frobenius
algebra with presentation $\mathfrak{p}=(\mathfrak{g}, \mathfrak{r})$.
Let $M$ be an irreducible $3$-manifold with $1$-manifold $\alpha \subset \partial M$. Then we have the isomorphism:
$$C(M,\alpha )\cong R\mathcal{I}(M,\alpha ,\mathfrak{g})/(\mathfrak{T}(M,\alpha ,\mathfrak{g}) +\mathfrak{C}(M, \alpha , \mathfrak{p}))=:\mathcal{D}(M, \alpha ).$$}

Note that in its definition $\mathcal{D}(M, \alpha )$ requires the presentation $\mathfrak{p}$. We omit the presentation from the notation here because the result shows that the quotient actually does not depend on it.

\vskip .1in

\noindent \textbf{Example 9.3.} If $\mathcal{F}$ is free with basis 
$\mathfrak{g}$ then  
$\mathfrak{C}(M, \alpha ,\mathfrak{p})=0$. Then in presentation 
9.2 all relations are tubing relations between 
incompressible marked surfaces.

\vskip .1in

We now define the isomorphism of the theorem and the tubing submodule. 
First consider $S\in \mathcal{S}(M,\alpha ,\mathfrak{g})$ 
and use 9.1. to expand a representative surface, also denoted $S$, 
as $R$-linear combination of incompressible 
marked surfaces. We do this in two steps. 

\vskip .1in

\noindent \textbf{Definition 9.4.} Let $M$ be a $3$-manifold with $\alpha \subset \partial M$ a closed $1$-manifold. Let $\mathfrak{a}\subset A$.
Then an $\mathfrak{a}$-\textit{pattern} $\Gamma $ is a finite graph with (i) a decomposition of its vertex set into a disjoint union of two sets, the black and the white vertices, (ii) a coloring of the components of the graph by elements 
of 
$\mathfrak{a}$ and, (iii) a one-to-one correspondence with the black vertices and the components of an irreducible surface in $M$. The set of $\mathfrak{a}$-patterns is denoted 
$\mathcal{P}(M,\alpha ,\mathfrak{a})$. 
 
\vskip .1in

For $\mathfrak{g}$ a generating set we call a $\mathfrak{g}$-pattern
just a \textit{pattern}. Define the \textit{projection}:
$$\mathfrak{h}: \mathcal{P}(M, \alpha ,\mathfrak{a})\rightarrow \mathcal{I}(M, \alpha ,1)$$

by assigning to each pattern the incompressible surface 
defined by its black vertices with trivial coloring of its components.

\vskip .1in

Now each presentation $\mathfrak{p}=(\mathfrak{g}, \mathfrak{a})$ of the algebra $A$ of a Frobenius algebra defines a 
\textit{state sum} map
$$\mathcal{K}: \mathcal{P}(M,\alpha ,\mathfrak{g})\rightarrow R\mathcal{I}(M,\alpha ,\mathfrak{g})/\mathfrak{C}(M, \alpha ,\mathfrak{p}).$$
Fix a Frobenius algebra for $\mathcal{F}$ such that
$$\Delta (1)=\sum_{i=1}^ru_i\otimes v_i.$$
Then a state $\sigma $ on $\Gamma $ assigns to each an element of 
$\{1,\ldots ,r\}$. We calculate the \textit{state evaluation} 
$\mathfrak{e}(\sigma ) \in \mathcal{I}(M,\alpha ,\mathfrak{g})$ as follows.
Choose an orientation of the graph. If an edge $e$ is running from 
a vertex $v$ to a vertex $w$ and $\sigma (e)=i$  
then assign to $v$ the element $u_i$ and to $w$ the element 
$v_i$. Finally assign the color $a$ of a component to any of its 
vertices. Multiply all the elements assigned to the vertices.
Consider the incompressible surface determined by the black 
vertices. Use the relations of $A$ to expand the colors in terms 
of the generating system. 
This is only defined up adding elements of 
$\mathfrak{C}(M, \alpha , \mathfrak{p})$.
The result is 
an element in $R\mathcal{I}(M,\alpha ,\mathfrak{g})/\mathfrak{C}(M, \alpha , \mathfrak{p})$ where the 
underlying surface is always the same but the colors in $\mathfrak{g}$ can 
change. Finally each white vertex $v$ is colored by some element 
$b_v\in A$. Multiply the element of $R\mathcal{I}(M,\alpha ,\mathfrak{g})/\mathfrak{C}(M, \alpha ,\mathfrak{p})$ determined above by 
$$\prod_{v \  \textrm{white vertex}}\varepsilon (b_v)\in R$$    
The resulting element is $\mathfrak{e}(\sigma )$. The sum over all states 
$\sigma $ defines $\mathcal{K}(\Gamma )$.
While $\mathfrak{e}(\sigma )$ in general will depend on the orientation
the full sum does not because of cocommutativity of $\Delta $.

\vskip .1in

\noindent \textbf{Example 9.5.} Consider the Bar-Natan system
$\mathcal{F}_{BN}$ with Frobenius 
algebra defined from $\{(1,X),(X,1)\}$. Then there are only two states that assign the dots to the different endpoints of edges. In this case the state itself can be identified with an orientation of the edges of the graph.
 
\vskip .1in

Next we assign to $S\in \mathcal{S}(M, \alpha ,\mathfrak{g})$ elements 
of $\Gamma \in \mathcal{P}(M,\alpha ,\mathfrak{g})$ such that 
$\mathcal{K}(\Gamma )$ are expansions of
$S$. Note that given $S$ we can apply neck cut compressions until the 
resulting surface is weakly incompressible. Because 
the neck cuttings always increase the Euler characteristic by $2$ 
this process stops after finitely many steps with a weakly incompressible 
surface with the \textit{history} of the neck cuttings defining the edges 
of the corresponding graph. 

Next consider the obvious map:
$$\mathfrak{d}: \mathcal{P}(M, \alpha ,\mathfrak{g})\times \mathcal{P}(M, \alpha ,\mathfrak{g})\rightarrow R\mathcal{I}(M,\alpha ,\mathfrak{g})/\mathfrak{C}(M, \alpha , \mathfrak{p})$$
defined by 
$$(\Gamma_1,\Gamma_2)\mapsto \mathcal{K}(\Gamma_1)-\mathcal{K}(\Gamma_2).$$  
For each $S\in \mathcal{S}(M, \alpha ,1)$ consider \textit{all} possible 
$\Gamma \in \mathcal{P}(M, \alpha ,\mathfrak{g})$ assigned to $S$ using the 
construction above. 
Let 
$$\mathcal{D}(S)\subset \mathcal{P}(M, \alpha ,\mathfrak{g})\times \mathcal{P}(M, \alpha ,\mathfrak{g})$$
be the set defined by all possible pairs of elements constructed from $S$ in this way and by assigning arbitrary colors to the different components of $\Gamma_1$ and $\Gamma_2$ in an arbitrary way such that corresponding components are colored in the same way.

\vskip .1in

\noindent \textbf{Definition 9.6.} Let $\mathfrak{T}(M, \alpha ,\mathfrak{g})$ be the submodule of $\mathcal{I}(M, \alpha ,\mathfrak{g})$ which is generated by the union of all sets
$\mathcal{K}(\mathfrak{d}(\mathcal{D}(S))$ for all $S\in \mathcal{S}(M, \alpha ,1)$. 

\vskip .1in

By construction the map:
$$\mathfrak{t}: \mathcal{S}(M,\alpha ,\mathfrak{g})\rightarrow \mathcal{D}(M, \alpha ) $$
is well-defined. 

\vskip .1in

\noindent \textit{Proof of} \textbf{9.2.} 
Using 8.1 we identify $C(M, \alpha )$ and $C(M, \alpha ,\mathfrak{g})$.
It is not hard to check that the 
relations (1')-(3') from section 8 are in the kernel of the $R$-module map 
$\mathfrak{t}: R\mathcal{S}(M, \alpha ,\mathfrak{g})\rightarrow \mathcal{D}(M, \alpha )$. For example consider a neck cutting relation. Then because the 
$\mathfrak{t}$-expansions are (by construction) independent of the way in which neck cutting relations are applied we can choose the expansions compatibly.
It is here important that we only have to consider neck cuttings on essential curves. Also sphere relations will contribute a factor that also is picked up in the state sum. The relations (1') induced from $\mathfrak{r}$ map to 
$\mathfrak{C}(M, \alpha , \mathfrak{p} )$ under the 
$\mathfrak{t}$-expansions. 
Conversely, the inclusion $\mathcal{I}(M, \alpha , \mathfrak{g})\rightarrow \mathcal{S}(M, \alpha ,\mathfrak{g})$ defines the $R$-module map:
$$\mathfrak{s}: R\mathcal{I}(M, \alpha ,\mathfrak{g})\rightarrow C(M, \alpha ,\mathfrak{g}).$$
It is immediate from the definitions that 
$$\mathfrak{T}(M, \alpha ,\mathfrak{g})+\mathfrak{C}(M, \alpha ,\mathfrak{r})$$
is contained in the kernel and thus $\mathfrak{s}$ induces an $R$-module map
$$\mathcal{D}(M, \alpha )\rightarrow C(M, \alpha ,\mathfrak{p}),$$
which is inverse to $\mathfrak{t}$. $\blacksquare$

\vskip .1in

\noindent \textbf{Remarks 9.7.} (a) Irreducibility is important in the discussion of skein modules of surfaces in $3$-manifolds. This has first been pointed out in \cite{AF}.
The result above is not true in general. For example one can show that in the 
Bar-Natan module of $M=(S^1\times S^2)\sharp (S^1\times S^2)$ the connected sum 
sphere represents an element $S\in C(M)$ with $4S=0$.

\noindent (b) For each $S_1, S_2\in \mathcal{I}(M, \alpha ,1)$ such that
$[S_1]=[S_2]\in H_2(M, \alpha ;\mathbb{Z}_2)$ there exists $S\in \mathcal{S}(M, \alpha ,1)$ and corresponding $(\Gamma_1, \Gamma_2)\in \mathcal{D}(S)$ which satisfy 
$\mathfrak{h}(\Gamma_i)=S_i$
for $i=1,2$. This foillows from 7.4. In particular $\mathcal{D}(S)\neq \emptyset $. 
Thus two incompressible surfaces within the same homology class
give rise to relations
$$\sum_{i=1}^{n_1}r_i^1(S_1,w_i^1)=\sum_{j=1}^{n_2}r_j^2(S_2,w_j^2)$$
for $w_i^1: S_1\rightarrow \mathfrak{g}$ and 
$w_j^2: S_2\rightarrow \mathfrak{g}$ suitable $\mathfrak{g}$-colorings 
and $r_i^1, r_2^j\in R$ for $i=1,\ldots n_1$ and $j=1,\ldots ,n_2$.
  
\vskip .1in

\noindent \textbf{Example 9.8.} Consider our version of the Gad-Naot skein module as discussed in 8.3 (c) with $R=\mathbb{Z}$.  Let $M=S^3\setminus int(N(K))$ where $N(K)\subset S^3$ is a tubular neighborhood of a nontrivial knot $K\subset S^3$. and $int$ is its interior. Then the boundary torus $S:=\partial N(K)\subset M$ is an incompressible surface representing an element $T\in C(M)$. Let $g$ be the tunnel number of the knot $k$. Then after attaching $g$ tubes to $S$ the resulting surface bounds a handlebody in $M$. Thus 
by $g$ applications of the neck cutting relation, 8.3 (c) and 
$X^2=t$, we conclude
$$(S,2^gX^g)=(1+(-1)^{g})^{g+1}t^{\frac{g}{2}}.$$
Thus for $g$ even we get
$$2^gt^{\frac{g}{2}}S=2^{g+1}t^{\frac{g}{2}}$$
or
$$2^gt^{\frac{g}{2}}(S-2)=0,$$
and for $g$ odd we get
$$2^gt^{\frac{g-1}{2}}(S,X)=0.$$
Tunnel numbers of knots are highly nontrivial and interesting invariants 
of the classical topology of knots in $3$-manifolds, see e.\ g.\ \cite{Mo}.

\end{document}